\documentclass{amsart}
\usepackage{amsthm}
\usepackage{amsmath}
\usepackage{amssymb}
\newtheorem*{theorem}{Theorem}
\newtheorem{remark}{Remark}
\def\O{{\mathcal O}}
\def\Ac{{\mathcal A}}
\def\Bc{{\mathcal B}}
\def\Mc{{\mathcal M}}

\def\Jc{{\mathcal J}}

\def\Hom{{\mathcal H om}}
\def\banica{{B\u anic\u a}}

\def\qed{\hfill$\Box$\vskip10pt}

\begin{document}

\title{Gorenstein multiple structures on smooth algebraic varieties}
\author{Nicolae Manolache}
\address{Institute of Mathematics "Simion Stoilow"\\
of the Romanian Academy \\
P.O.Box 1-764
Bucharest, RO-014700
}
\email{nicolae.manolache@imar.ro}
\date{}
\begin{abstract}
We characterize the Gorenstein nilpotent scheme structures  on a smooth algebraic variety as support, 
in terms of a duality property of the graded objects associated to two canonical filtrations.
\end{abstract}
\maketitle

\section{Preliminaries}

We present here, for the sake of the reader, the frame of our note.

Let $X$ be a smooth connected algebraic variety over an algebraically closed field $k$.
A locally Cohen-Macaulay scheme $Y$ is called a \emph{multiple structure on $X$} if the subjacent 
reduced scheme $Y_{red}$ is $X$. In this case all the local rings of $Y$ have the same multiplicity
(cf. \cite{M1}, \cite{M2}), which is called \emph{the multiplicity of $Y$}. To a given multiple structure $Y$
on $X$ one associates canonically three filtrations. To set the frame of the following considerations, 
let $I$ be the (sheaf) ideal of $X$ in $Y$ and let $m$ be the positive integer such that $I^m \neq 0$,
$I^{m+1}=0$. The three filtrations are:\\
1. Let $I^{(\ell)}$ be the ideal obtained throwing away the embedded components of $I^\ell$ and let 
$Z_\ell$ be the corresponding scheme. This  gives the 
\emph{\banica -Forster filtration} (cf. \cite{BF2}):

\begin{displaymath}
\begin{array}{ccccccccccccc}
 \O_Y=I^{(0)} & \supset & I= I^{(1)} & \supset &  I^{(2)} & \supset &  \ldots &  \supset &   I^{(m)} & \supset & I^{(m+1)}=0\medskip \\
&  & X=Z_1 & \subset & Z_2 & \subset &  \ldots &  \subset & Z_m & \subset & Z_{m+1} =  Y
\end{array}
\end{displaymath}
$Z_\ell$ are not, in general, Cohen-Macaulay. But this is true if ${\rm dim}(X)=1$.
The graded associated object $\Bc (Y)=\bigoplus_{\ell =0}^m I^{(\ell)} / I^{(\ell +1)} $ is naturally a 
graded $\O_X$-algebra. If the schemes $Z_\ell$ are Cohen-Macaulay, the graded components of $\Bc (Y)$ are locally 
free sheaves on $X$.\\
2. Let $X_\ell $ be defined by $I_\ell =0:I^{m+1-\ell}$. Again, if ${\rm dim}(X)=1$, $X_\ell $ are Cohen-Macaulay.
This is also true if $Y$ is locally complete intersection of multiplicity
at most $6$ (cf. \cite{M2}).
In general this is not always the case. When $X_\ell$ are Cohen-Macaulay, the quotients $I_\ell /I_{\ell +1}$ are locally free sheaves on $X$. This filtration was considered in \cite{M1}.\\
3. Let $Y_\ell$ be the scheme given by $J_\ell=0:I_{m+1-\ell}=0:(0:I^\ell)$. When $X_\ell$ is Cohen-Macaulay, $Y_\ell$
has the same property. The graded object $\Ac (Y)=\bigoplus_{\ell =0}^m \Jc_\ell / \Jc_{\ell +1}$ is a graded 
$\O _X$-algebra and $\Mc (Y) =\bigoplus_{\ell =0}^m I_\ell / I_{\ell+1}$ is a graded $\Ac(Y)$-module.
This filtration was considered in \cite{M2} \\
The system of the graded components ($\Ac_0(Y),\ldots  \Ac_m(Y);\Mc_0(Y),\ldots  \Mc_m(Y))$ is called 
\emph{the type of $Y$}. $Y$ is called \emph{of free type} when all the graded pieces are locally free.
As already remarked, in dimension $1$, or if $Y$ is lci of multiplicity up to $6$, this is  the case.

Recall some properties:\\
1) In general the above filtration are different.
Take for instance $X=Spec(k)$, $Y=Spec(k[x,y]/(x^3,xy,y^4))$\\
2) $Z_\ell \subset Y_\ell \subset X_\ell$ \\
2') there are canonical morphisms: $\Bc(Y)\to \Ac(Y)\to \Mc (Y)$\\
3) The multiplications
\begin{displaymath}
\begin{array}{lll}
\Ac_{\ell_1}\otimes \Ac_{\ell_2} & \to & \Ac_{\ell_1+\ell_2}\medskip\\
\Ac_{\ell_1}\otimes \Mc_{\ell_2} & \to & \Mc_{\ell_1+\ell_2}
\end{array}
\end{displaymath}
are never the zero maps for $\ell_1, \ell_1 \ge 0$, $\ell_1+\ell_2 \le m$ (cf. \cite{M2}.\\
4) One has the exact sequences:
\begin{displaymath}
\begin{array}{c}
0\to \Mc_\ell(Y)\to \O_{X_{\ell+1}} \to \O_{X_\ell} \to 0\medskip\\
0\to \Ac_\ell(Y)\to \O_{Y_{\ell+1}} \to \O_{Y_\ell} \to 0
\end{array}
\end{displaymath}
5) If $Y$ is Gorenstein of free type, then $X_\ell $ and $Y_{m+1-\ell}$ are locally algebraically linked
(cf. \cite{M1}). In particular one has the exact sequences:
\begin{displaymath}
\begin{array}{c}
0\to \omega_{X_{m+1-\ell}}\otimes \omega_Y^{-1}\to \O_Y \to \O_{Y_\ell} \to 0\medskip\\
0\to \omega_{Y_{m+1-\ell}}\otimes \omega_Y^{-1}\to \O_Y \to \O_{X_\ell} \to 0
\end{array}
\end{displaymath}
6) If $Y$ is Gorenstein of free type, then (cf \cite{M2}):\\
\begin{displaymath}
\begin{array}{l}
(a)\  \rm{rank}\ \Ac_\ell(Y)=\rm{rank}\ \Mc_{m-\ell}(Y)\medskip\\
(b)\  \Ac_\ell(Y)=\Mc _\ell (Y) \rm{\ iff\ } \rm{rank}\ \Ac_\ell(Y)=\rm{rank}\ \Ac_{m-\ell}(Y)
\end{array}
\end{displaymath}

\section{Main theorem}
The aim of this theorem is to "explain" the equality 6)(a) from above.
In fact one gives a characterization of the Gorenstein multiple structures of free type
on a smooth support which generalizes the result from \cite{Bo}.
\begin{theorem}
Let $Y$ be a free type Cohen-Macaulay multiple structure on a smooth support $X$.\\
Then $Y$ is Gorenstein if and only if the following conditions are fulfilled:\\
(a) $\Ac_m$ and $\Mc_m$ are line bundles\\
(b) $\Ac_m=\Mc_m$  \\
(c) The canonical maps:
\begin{displaymath}
\Ac_\ell \to \Hom _{\O_X}(\Mc_{m-\ell},\Mc_m)\cong \Mc_{m-\ell}^\vee\otimes \Mc_\ell
\end{displaymath}
are isomorphisms.
\end{theorem}
\emph{Proof.}
Suppose $Y$ is Gorenstein. Dualizing (i.e. applying $\Hom (?,\omega _Y)$) the exact sequences:
\begin{displaymath}
\begin{array}{c}
0\to \Ac_m(Y)\to \O_Y \to \O_{Y_m} \to 0\medskip\\
0\to \Mc_m(Y)\to \O_Y \to \O_{X_m} \to 0
\end{array}
\end{displaymath}
and taking the restrictions to $X$ one gets:
\begin{displaymath}
\omega_Y\lvert _X\cong \Ac_m(Y)^\vee\otimes \omega_X\cong \Mc_m(Y)^\vee \otimes \omega_X
\end{displaymath}
and so (a) and (b) are fulfilled.\\
Dualizing the exact sequence:
\begin{displaymath}
0\to \Mc_{m-\ell}(Y)\to \O_{X _{m-\ell+1}}\to \O_{X_{m-\ell}} \to 0
\end{displaymath}
one gets the exact sequence:
{\scriptsize
\begin{displaymath}
\begin{array}{cc}
0\to \Hom_{\O_Y}(\O_{X_{m-\ell}},\omega_Y)\to \Hom_{\O_Y}(\O_{X_{m-\ell+1}},\omega_Y)\to &
\Hom_{\O_Y}(\Mc_{m-\ell}(Y),\omega_Y)\to 0\medskip\\
& \Vert \medskip \\
& \Hom_{\O_X}(\Mc_{m-\ell}(Y),\omega_X)
\end{array}
\end{displaymath}
}
which tensored with $\omega _Y^{-1}$ gives the exact sequence:

{\tiny
\begin{displaymath}
\begin{array}{ccc}
0\to \Hom_{\O_Y}(\O_{X_{m-\ell}},\O_Y)\to & \Hom_{\O_Y}(\O_{X_{m-\ell+1}},\O_Y)\to &
\Hom_{\O_Y}(\Mc_{m-\ell}(Y),\omega_Y)\otimes \omega_Y^{-1}\to 0\medskip\\
\Vert & \Vert & \Vert \medskip \\
0\to 0:(0:I^{\ell+1}) & 0:(0:I^\ell)  & \Hom_{\O_X}(\Mc_{m-\ell}(Y),\omega_X\otimes \omega _Y^{-1}) \to 
0\medskip\\
& & \Vert \\
& & \Hom_{\O_X}(\Mc_{m-\ell}(Y),\Mc_m(Y))
\end{array}
\end{displaymath}
}
and so (c) is fulfilled.

Assume now (a), (b), (c) fulfilled. The Gorenstein property being local and taking the
completions of the local rings, we are rduced to the following situation:

$B$ is a local complete ring, $B_{\rm red}=A\cong k[[X_1,\ldots ,X_n]]$ and the surjection $p:B\to
A$ admits a retract
$i:A\to B$ which makes $B$ an $A$-algebra. If $A_\ell:=B/(0:I^{m+1-\ell})$, $B_\ell:= B/(0:(0:I^\ell))$, 
one has the split exact sequences:
\begin{displaymath}
\begin{array}{ccc}
0\to \Ac_m &\to B_{m+1} &\to B_m \to 0\medskip \\
0\to \Ac_{m-1} &\to B_m &\to B_{m-1} \to 0\medskip\\
\vdots & \vdots & \vdots \medskip \\
0\to \Ac_1 &\to B_2 &\to B_1 \to 0
\end{array}
\end{displaymath}
Observe that $B_1=B/(0:(0:I))=A=\Ac_0$; one obtains 
$B\cong \bigoplus_{\ell=0}^M\Ac_\ell$ and this is an isomorphism of $A$-algebras.

Similarly, one has a morphism of $A$-modules $B\cong \bigoplus_{\ell=0}^m\Mc_\ell$.

Then:
\begin{displaymath}
\begin{array}{l}
\omega_B\cong Hom_A(B,A)\cong Hom_A(\oplus_{\ell=0}^m M_\ell, A) \cong Hom_A(\oplus_{\ell=0}^m M_\ell, M_m) \cong\medskip \\
\oplus_{\ell=0}^m Hom_A(M_{m-\ell}, M_m)\cong \oplus B_\ell \cong B \ .
\end{array}
\end{displaymath}
This shows that $B$ is Gorenstein.
\qed
\begin{remark}
The case of embedded multiple structures, (as in \cite{BF2}, \cite{M1}, \cite{M2}) i.e. $X\subset Y \subset P$, with $X$, $P$ 
smooth connected algebraic varieties and $Y$ a multiple structure on $X$, leads to similar filtrations. Denote now by $I$ 
the ideal of $X$ in $P$, by $J$ the ideal of $Y$ in $P$ and suppose $I^{m+1}\subset J$, $I^m \not\subset J$.
Then the filtrations are (we consider now only the ideals):\\
1. Let $I^{(k)}$ be the ideal obtained throwing away the embedded components of $I^{(k)}+J$\\
2. Let $I_\ell = J:I^{m+1-\ell}$\\
3. Let $J_\ell =J:(J:I^\ell)$

All the above considerations applies also in this case. 
\end{remark}
\begin{remark}
In the case of quasiprimitive structures, i.e. when ${\rm rank }\Ac_\ell = {\rm rank }\Mc_\ell =1$ for all $\ell$,
all the above filtrations are equal, whence $\Bc=\Ac=\Mc$, so one can express the conditions in the theorem only in terms of \banica - Forster filtration. This was done in \cite{Bo} for the case of quasiprimitive multiple structures on smooth curves in a threefold.
\end{remark}
\begin{remark}
The above duality gives a direct explanation of the identities in the Chern classes of the bundles which appear in various constructions in \cite{M2} (e.g. 4.14, 4.16. loc.cit.).
\end{remark}
{\bf Aknowledgment} The author was partially supported by the Humboldt Foundation
and the CNCSIS 33079/2004 contract.

\noindent Nicolae Manolache\\
Institute of Mathematics "Simion Stoilow"\\
of the Romanian Academy \\
P.O.Box 1-764
Bucharest, RO-014700

\noindent e-mail: nicolae.manolache@imar.ro
\end{document}